
\baselineskip=14pt
\parskip=10pt

\magnification=\magstephalf

\def\1{{\overline{1}}}
\def\2{{\overline{2}}}
\parindent=0pt
\overfullrule=0in

\def\frac#1#2{{#1 \over #2}}
\centerline
{\bf 
Linear-Time and Constant-Space Algorithms to compute Multi-Sequences  
}
\centerline
{\bf 
that arise in Enumerative Combinatorics (and Elsewhere)
}

\bigskip
\centerline
{\it Shalosh B. EKHAD and Doron ZEILBERGER}
\bigskip

{\bf Abstract}: How many ways, {\bf exactly},  can a Chess King, always moving forward, i.e. with the set of steps $\{[1,0],[0,1],[1,1]\}$, walk from the origin, $(0,0)$, to the lattice point $(100000,200000)$?
Thanks to the amazing Apagodu-Zeilberger extension of the Almkvist-Zeilberger algorithm, adapted in this article for combinatorial applications,
this  $104492$-digit number, viewable from  \hfill\break
{\tt http://www.math.rutgers.edu/\~{}zeilberg/tokhniot/oPureRecRat2.txt},
can be computed in less than 33 seconds. But not just this particular number! Many other numbers that come up in enumerative combinatorics, and elsewhere, can be computed just as efficiently.

{\bf Maple Packages}

This article is accompanied by two  Maple packages, {\tt PureRec.txt}, and {\tt PureRecRat.txt}, available from 

{\tt https://sites.math.rutgers.edu/\~{}zeilberg/tokhniot/PureRec.txt} \quad  and

{\tt https://sites.math.rutgers.edu/\~{}zeilberg/tokhniot/PureRecRat.txt} \quad  .

The web-page of this article,

{\tt https://sites.math.rutgers.edu/\~{}zeilberg/mamarim/mamarimhtml/pure.html} \quad ,

contains input and output files, referred to in this paper.

The purpose of this short article is to be a {\bf pointer} to the much longer Maple packages (written by DZ) and numerous output files (generated by SBE) given in the above mentioned
web-page. Recall that the Almkvist-Zeilberger algorithm [AlZ] (see [D] for an engaging exposition)
produces a {\bf pure} {\it linear} recurrence equation with {\it polynomial coefficients}, satisfied by the Taylor coefficients of many kinds of generating functions. All we need is that
$R'(x)/R(x)$ is a rational function of $x$. This includes all rational functions, their powers, $R(x)^{\alpha}$ (either for symbolic or numeric exponent $\alpha$)
and all functions of the form $e^{P(x)}$ where $P(x)$ is a polynomial, and many other functions besides.
This follows from the fact, that thanks to Cauchy, it can be written as a contour integral $\frac{1}{2 \pi i}\,\int_{|x|=1} \frac{R(x)}{x^{n+1}} \, dx$. In fact, it can do much more, it can handle
any expressions of the form   $\int_{|x|=1} R(x,n) \, dx$, where $R(x,n)$ satisfies the condition that both $\frac{\frac{d}{dx}R(n,x)}{R(n,x)}$ and $\frac{R(n+1,x)}{R(n,x)}$ are {\bf rational functions} of
the discrete variable $n$ and the continuous variable $x$.

Since the recurrence is {\bf finite} order (often fairly low), this enables a very fast linear-time and constant-memory computation of such sequences. If the order of the recurrence is $L$, one only has to `remember' the last $L$ values.
Of course, since the numbers (usually) grow exponentially, the {\it bit-size} memory is linear, but still it is very good!

After dividing by the leading term, one can write the recurrence for the sequence, let's call is $a(n)$ as
$$
a(n)= \sum_{i=1}^{L} f_i(n) a(n-i) \quad,
$$
where $f_i(n)$ are {\bf rational functions} of $n$. There is a potential {\it issue} of `dividing by zero', but if $n_0$ is the largest positive integer root of the denominators, one can directly compute the first $n_0$ values and start
the recurrence there.

The  amazing {\bf Multi-variable Almkvist-Zeilberger algorithm}, designed by Moa Apagodu and Doron Zeilberger [ApZ]  extends this to several variables. If $a(n_1, \dots, n_d)$ is such a multi-sequence, gotten, for example.
from the Taylor coefficients of a multi-variable rational function in $d$ variables (and many other kinds of  multi-sequences) , then the algorithm finds $d$ {\bf pure} recurrences, one for each discrete axis:

$$
a(n_1, \dots, n_d)= \sum_{i=1}^{L_j} f^{(j)}_i(n) a(n_1, \dots, n_{j-1}, n_{j} -i, n_{j+1}, \dots, n_d ) \quad,
$$
where $f^{(j)}_i(n_1, \dots, n_d)$ are {\bf rational functions} of $n_1, \dots, n_d$. 
Then to compute a specific value,  $a(n_1, \dots, n_d)$, one constructs a discrete-path from the origin $[0,\dots,0]$ to the lattice point $(a_1, \dots, a_d)$ 
with $d$ segments, one for each direction. One such path does not suffice, since in order to say, climb-up in the last axis $n_d$, we need $L_d$ initial values, so 
we need recursively to take care of them $L_d$ `initial conditions', each requiring its own path in $(d-1)$-discrete space,
but one only needs finitely many paths, making it linear time and constant memory.

Alas, now the {\bf division by zero} is not so easy to overcome. We chose to only use the most obvious generic paths (see the Maple code for procedure {\tt EvalScheme} in the Maple package \hfill\break
{\tt https://sites.math.rutgers.edu/\~{}zeilberg/tokhniot/PureRec.txt} \quad), and our program returns FAIL if one encounters a division by $0$. However, we believe that this can be overcome by finding
other paths that stay-away from dividing by $0$. 

{\bf Sample Output}

For five random examples of sequences that are Taylor coefficients of one-variable  functions of the form $1/P(x)^{\frac{1}{3}}$, where $P(x)$ is a randomly generated polynomial, see:

{\tt https://sites.math.rutgers.edu/\~{}zeilberg/tokhniot/oPureRec1.txt} \quad .

For five random examples of  sequences that are Taylor coefficients of one-variable functions that are exponentials of polynomials, see:

{\tt https://sites.math.rutgers.edu/\~{}zeilberg/tokhniot/oPureRec1a.txt} \quad .

For five random examples  of  bi-sequences that are Taylor coefficients of two-variable  functions of the form $1/P(x_1,x_2)^{\frac{1}{3}}$, where $P(x_1,x_2)$ is a randomly generated polynomial of $x_1,x_2$, see:

{\tt https://sites.math.rutgers.edu/\~{}zeilberg/tokhniot/oPureRec2.txt} \quad .

For five random examples  of  bi-sequences that are Taylor coefficients of
functions  that are exponentials of polynomials of two variables, see

{\tt https://sites.math.rutgers.edu/\~{}zeilberg/tokhniot/oPureRec2a.txt} \quad .

For three random examples of triple-sequences that are Taylor coefficients of a
three-variable function  of the form $1/P(x_1,x_2,x_3)^{\frac{1}{3}}$, where $P(x_1,x_2,x_3)$ is a randomly generated polynomial of $x_1,x_2,x_3$, see:

{\tt https://sites.math.rutgers.edu/\~{}zeilberg/tokhniot/oPureRec3.txt} \quad .

Now it takes much longer to {\it generate} the scheme, but once found, it is still very fast to compute specific values, even those very far from the origin,

For three random examples of tri-sequences that are Taylor coefficients of functions that are exponentials of polynomials of three variables, see

{\tt https://sites.math.rutgers.edu/\~{}zeilberg/tokhniot/oPureRec3a.txt} \quad . 

Here we encountered some singularities, so our (incomplete) implementation sometimes (but fairly rarely) returns FAIL.

{\bf Applications to the Enumeration of Lattice Path in the Plane with Many choices of Atomic Steps}

Let {\tt St} be an arbitrary set of positive steps, e.g. for the forward-moving Chess King it is $\{[1,0],[0,1],[1,1]\}$. How can we compute fast the number of ways
of walking from the origin to $[a,b]$ using the members of {\tt St}?

The lattice-path oriented companion Maple package

{\tt https://sites.math.rutgers.edu/\~{}zeilberg/tokhniot/PureRecRat.txt} \quad,

has a procedure

{\tt Walk2Dpaper(St,K)} \quad,

that inputs an arbitrary set of positive steps in 2-dimensions, and a large positive integer K (for illustrative purposes) and outputs
the scheme for the bi-variate function $F(a,b)$, defined as the number of walks from $[0,0]$ to $[a,b]$ using the steps in {\tt St}, and if
no singularity is encountered, outputs $F(K,2K)$. It also outputs the one-variable linear recurrence for the diagonal sequence
$\{F(a,a)\}$ and uses it to compute $F(2K,2K)$ (now there are never singularities). For the sake of Neil Sloane, it also spits
out the first $30$ terms of the diagonal sequence.

You are welcome to try out procedure {\tt Walk2Dpaper(St,K)} with $St=\{[0,1],[1,0],[1,1]\}$ and $K=100000$. getting the $104922$-digit number mentioned in the abstract.

For $29$ such theorems see:

{\tt https://sites.math.rutgers.edu/\~{}zeilberg/tokhniot/oPureRecRat1.txt} \quad  .

{\bf Christoph Koutschan's even better program (alas, in  Mathematica)}

A Mathematica implementation of the multi-variate Almkvist-Zeilberger algorithm is contained, {\it inter alia}, in Christoph Koutschan's very
efficient and versatile Mathematica program [K]. It would be nice if he would use the simple ideas in our paper to interface them
and create such a {\it multi-sequence calculator}.

{\bf References}

[AlZ]  Gert Almkvist and Doron Zeilberger, {\it The method of differentiating Under The integral sign}, J. Symbolic Computation {\bf 10} (1990), 571-591. \hfill \break
{\tt https://sites.math.rutgers.edu/\~{}zeilberg/mamarim/mamarimPDF/duis.pdf} \quad .

[ApZ] Moa Apagodu and Doron Zeilberger,
{\it Multi-Variable Zeilberger and Almkvist-Zeilberger Algorithms and the Sharpening of Wilf-Zeilberger Theory},
Adv. Appl. Math. {\bf 37} (2006), 139-152. [Special issue in honor of Amitai Regev] \hfill\break
{\tt https://sites.math.rutgers.edu/\~{}zeilberg/mamarim/mamarimhtml/multiZ.html} \quad .

[D] Robert Dougherty-Bliss, {\it Integral recurrences from A to Z}, American Mathematical Monthly, to appear. \hfill\break
{\tt https://arxiv.org/abs/2102.10170} \quad .

[K] Christoph Koutschan, {\it   Advanced applications of the holonomic systems approach}, 
PhD thesis, Research Institute for Symbolic Computation (RISC), Johannes Kepler University, Linz, Austria, 2009.\hfill\break
{\tt http://www.koutschan.de/publ/Koutschan09/thesisKoutschan.pdf}, \hfill\break

\bigskip
\hrule
\bigskip
Shalosh B. Ekhad and Doron Zeilberger, Department of Mathematics, Rutgers University (New Brunswick), Hill Center-Busch Campus, 110 Frelinghuysen
Rd., Piscataway, NJ 08854-8019, USA. \hfill\break
Email: {\tt [ShaloshBEkhad, DoronZeil] at gmail dot com}   \quad .

{\bf Exclusively published in the Personal Journal of Shalosh B. Ekhad and Doron Zeilberger and arxiv.org}

{\bf March 9, 2022}

\end